\documentclass[a4paper,reqno]{amsart}
\usepackage{amssymb, amsmath, amscd}
\usepackage{color}
\def\HLP{[0,\infty)}\def\HLN{(-\infty,0]}

\def\black{\color{black}}
\newtheorem{theorem}{Theorem}[section]

\newtheorem{example}[theorem]{Example}
\makeatletter
\renewcommand{\theequation}{%
\thesection.\alph{equation}} \@addtoreset{equation}{section}
\makeatother
\begin{document}
\makeatletter
\renewcommand{\theequation}{%
\thesection.\alph{equation}} \@addtoreset{equation}{section}
\makeatother
\title[Proper curves]
{Curves in $\mathbb{R}^n$ with finite total first curvature arising
from the solutions of an ODE}
\author
{P. Gilkey, C. Y. Kim, H. Matsuda, J. H. Park, and S. Yorozu}
\address{PG: Institute of Theoretical Science, University of Oregon,
Eugene OR 97403 USA}
\email{gilkey@uoregon.edu}
\address{CYK \& JHP:Department of Mathematics,
Sungkyunkwan University,
Suwon, 440-746, Korea}
\email{intcomplex@skku.edu, parkj@skku.edu}
\address{HM: Department of Mathematics,
Kanazawa Medical University, Uchinada, Ishikawa, 920-0293, Japan}
\email{matsuda@kanazawa-med.ac.jp}\address{SY: Miyagi University of
Education, Sendai, Miyagi, 980-0845, Japan}
\email{s-yoro@staff.miyakyo-u.ac.jp}
\begin{abstract} We use the solution set
of a real ordinary differential equation which has
order $n\ge2$ to construct a smooth curve $\sigma$ in $\mathbb{R}^n$.
We describe when $\sigma$ is a proper embedding of infinite length
with finite total first curvature.
\end{abstract}
\subjclass[2010]{53A04}
\keywords{finite total curvature, ordinary differential equation,
proper embedded curve}
\maketitle

\section{Statement of results}
Let $\sigma$ be an immersion of $\mathbb{R}$ into $\mathbb{R}^n$
for $n\ge2$; if $\sigma$ is proper, then
the length of $\sigma$ is then necessarily infinite.
The first curvature $\kappa$ and the total first curvature $\kappa[\sigma]$
are given, respectively, by:
\begin{equation}\label{eqn-1.a}
\kappa:=\frac{||\dot\sigma\wedge\ddot\sigma||}{||\dot\sigma||^3}\text{ and }
\kappa[\sigma]:=\int_{-\infty}^\infty\kappa||\dot\sigma||dx\,.
\end{equation}

\subsection{History}
Fenchel \cite{F40} showed that a simple closed curve in $\mathbb{R}^3$
had $\kappa[\sigma]\ge2\pi$.
F\'ary \cite{F49} and Milnor \cite{M53} showed the total curvature of
any knot is greater than $4\pi$.
Castrill\'on L\'opez and Fern\'andez Mateos \cite{CM10},
and Kondo and Tanaka \cite{KT11} have examined the global properties
of the total curvature of a curve in an arbitrary Riemannian manifold.
The total curvature of open plane curves of fixed length in $\mathbb{R}^2$
was studied by Enomoto \cite{E00}. The analogous question
for $S^2$ was examined by Enomoto and Itoh \cite{EI04, EI04a}.
Enomoto, Itoh, and Sinclair \cite{EIS08} examined curves in $\mathbb{R}^3$.
We also refer to related work of Sullivan \cite{S08}.
Borisenko and Tenenblat \cite{BT12} studied the problem in
Minkowski space.
Buck and Simon \cite{BS07} and Diao and Ernst \cite{DE07}
studied this invariant
in the context of knot theory,
and Ekholm \cite{E06} used this invariant
in the context of algebraic topology. Alexander, Bishop, and Ghrist \cite{ABG10}
extended these notions to more general spaces than smooth manifolds.
The total curvature also appears in the study of Plateau's
problem -- see the discussion in Desideri and Jakob \cite{DJ13}.

The literature on the subject is a vast one and we have only cited a few representative
papers to give a flavor for the subject.

The papers cited above focused
on closed curves, polygonal curves, knotted curves, curves with
fixed endpoints, curves with finite length and pursuit curves. This present paper
deals, by contrast, with properly embedded curves which arise from
ordinary differential equations.
We take as our starting point a real constant coefficient ordinary differential
operator $P$ of degree $n\ge2$ of the form:
$$P(y):=y^{(n)}+c_{n-1}y^{(n-1)}+\dots+c_0y\,.$$
Let $\mathcal{S}=\mathcal{S}_P$ be the solution set,
let $\mathcal{P}=\mathcal{P}_P$ be the associated characteristic polynomial,
and let $\mathcal{R}=\mathcal{R}_P$ be the roots of $\mathcal{P}$:
\begin{eqnarray*}
&&\mathcal{S}:=\{y:P(y)=0\},\\
&&\mathcal{P}(\lambda):=
\lambda^n+c_{n-1}\lambda^{n-1}+\dots+c_0,\\
&&\mathcal{R}:=
\{\lambda\in\mathbb{C}:\mathcal{P}(\lambda)=0\}\,.
\end{eqnarray*}
We suppose for the moment that all the roots of $\mathcal{P}$ have
multiplicity $1$ and enumerate the roots of $\mathcal{P}$ in the form:
$$\mathcal{R}=\{s_1,\dots,s_k,t_1,\bar t_1,\dots,t_\ell,\bar t_\ell\}
\text{ for }k+2\ell=n$$
where $s_i\in\mathbb{R}$ for $1\le i\le k$
and where $t_j=a_j+\sqrt{-1}b_j$ with $b_j>0$ for $1\le j\le\ell$.
Since we have assumed that all the roots are distinct, the standard basis for
$\mathcal{S}$ is given by the functions
\begin{equation}\label{eqn-1.b}
\begin{array}{lrr}
\phi_1:=e^{s_1x},&\dots,&\phi_k:=e^{s_kx},\\
\phi_{k+1}:=e^{a_1x}\cos(b_1x),&
\phi_{k+2}:=e^{a_1x}\sin(b_1x),&\dots,\\
\phi_{n-1}:=e^{a_\ell x}\cos(b_\ell x),&\phi_n:=e^{a_\ell x}\sin(b_\ell x)\,.
\end{array}\end{equation}
Of course, if all the roots are real, then $k=n$ and we omit the functions
involving $\cos(\cdot)$ and $\sin(\cdot)$.
Similarly, if all the roots are complex, then $k=0$ and we omit the
pure exponential functions.
We define the associated curve $\sigma_P:\mathbb{R}\rightarrow\mathbb{R}^n$
by setting:
$$\sigma_P(x):=(\phi_1(x),\dots,\phi_n(x))\,.$$

\subsection{The length of the curve $\sigma_P$}
Let $\Re(\lambda)$ denote the real part of a complex number $\lambda$. Define:
\begin{eqnarray*}
&&r_+(P):=\max_{\lambda\in\mathcal{R}}\Re(\lambda)
=\max(s_1,\dots,s_k,a_1,\dots,a_\ell),\\
&&r_-(P):=\min_{\lambda\in\mathcal{R}}\Re(\lambda)
=\min(s_1,\dots,s_k,a_1,\dots,a_\ell)\,.
\end{eqnarray*}
The numbers $r_\pm(P)$ control the growth of $||\sigma_P||$ as
$x\rightarrow\pm\infty$. Section~\ref{sect-2} is devoted to the proof of the following result:

\begin{theorem}\label{thm-1.1}
 Assume that all the roots of $\mathcal{P}$ are simple. If $r_+(P)>0$, then
$\sigma_P$ is a proper embedding of $\HLP$ into $\mathbb{R}^n$ with
infinite length.
If $r_-(P)<0$, then
$\sigma_P$ is a proper embedding of $\HLN$ into $\mathbb{R}^n$ with
infinite length.
\end{theorem}

\subsection{The total first curvature}
If $\sigma$ is a curve, then the first curvature $\kappa$ of
Equation~(\ref{eqn-1.a}) is a local invariant of the curve which
does not depend on the parametrization.
Let $\rho(x)$ be the radius of the best circle approximating $\sigma$ at $x$.
We then have that
$\kappa=\rho^{-1}$. If we were to reparametrize $\sigma$ by
arc length, then $\ddot\sigma$ would be perpendicular to $\dot\sigma$. In
particular in that setting, $\kappa$ vanishes identically if and only if
$\ddot\sigma=0$ or equivalently if $\sigma$ is a geodesic. Thus $\kappa[\sigma]$
is also called the {\it geodesic curvature} and is a very important invariant
of the curve.

Let $\kappa_+$ (resp. $\kappa_-$) be the total first curvature of $\sigma_P$ on
$\HLP$ (resp. on $\HLN$):
$$
\kappa_+:=
\int_0^\infty\frac{||\dot\sigma_P\wedge\ddot\sigma_P||}{||\dot\sigma_P||^2}dx
\text{ and }
\kappa_-:=
\int_{-\infty}^0\frac{||\dot\sigma_P\wedge\ddot\sigma_P||}{||\dot\sigma_P||^2}dx\,.
$$
We normalize the labeling of the roots so that
$$
s_1>s_2>\dots>s_k\text{ and }
a_1\ge \dots \ge a_\ell\,.
$$
We then have $r_+=\max(s_1,a_1)$ and
$r_-=\min(s_k,a_\ell)$.
Section~\ref{sect-3} is devoted to the proof of the following result:
\begin{theorem}\label{thm-1.2}
Assume that all the roots of $\mathcal{P}$ are simple,
that $r_+(P)>0$, and that $r_-(P)<0$.
\begin{enumerate}
\item If $s_1>a_1$, then $\kappa_+<\infty$;
otherwise, $\kappa_+=\infty$.
\item If $s_k<a_\ell$, then $\kappa_-<\infty$;
otherwise $\kappa_-=\infty$.
\end{enumerate}
\end{theorem}

We note that if there are no complex roots, then $s_1>0$ and $s_k<0$ and
we may conclude that
$\kappa_+$ and $\kappa_-$ are finite. This is quite striking as these curves
are, obviously, not straight lines. On the other hand, if there are no real roots,
then $a_1>0$ and $a_\ell<0$ and we may conclude that
$\kappa_+$ and $\kappa_-$ are infinite.
\subsection{Examples} Section~\ref{sect-4} treats several families
of examples. We construct
examples where $\kappa_+$ and $\kappa_-$ are both finite, where $\kappa_+$ is
finite but $\kappa_-$ is infinite, where $\kappa_+$ is infinite but $\kappa_-$ is finite,
and where both $\kappa_+$ and $\kappa_-$ are infinite.

\subsection{Changing the basis}
We have considered the standard basis for $\mathcal{S}$ to define the curve
$\sigma_P$. More generally, let $\Psi:=\{\psi_1,\dots,\psi_n\}$ be an arbitrary basis
for $\mathcal{S}$. We define:
$$\sigma_{\Psi,P}(x):=(\psi_1(x),\dots,\psi_n(x))\,.$$
In Section~\ref{sect-5}, we extend Theorem~\ref{thm-1.1} and
Theorem~\ref{thm-1.2} to this setting and verify that the properties
we have been discussing are properties of the solution space
$\mathcal{S}$ and not of the particular basis chosen:

\begin{theorem}\label{thm-1.3}
Assume that all the roots of $\mathcal{P}$ are simple,
that $r_+(P)>0$, and that $r_-(P)<0$.
Then $\sigma_{\Psi,P}$ is a proper embedding of $\HLP$ and of $\HLN$
into $\mathbb{R}^n$ with
infinite length.
\begin{enumerate}
 \item If $s_1>a_1$, then $\kappa_+[\sigma_{\Psi,P}]<\infty$;
otherwise, $\kappa_+[\sigma_{\Psi,P}]=\infty$.
\item If $s_k<a_\ell$, then $\kappa_-[\sigma_{\Psi,P}]<\infty$;
otherwise $\kappa_-[\sigma_{\Psi,P}]=\infty$.
\end{enumerate}\end{theorem}

\subsection{Roots with multiplicity greater than 1}
If we have roots of multiplicity
greater than $1$, then powers of $x$ arise. For example, if we consider the equation
$y^{(n)}=0$, then
$$\mathcal{S}=\operatorname{Span}\{1,\ x,\ \dots,\ x^{n-1}\}\,.$$
More generally, if $s$ is a real eigenvalue
of multiplicity $\nu\ge2$, then we must consider the family of functions:
\begin{equation}\label{eqn-1.c}
\{\phi_{s,0}:=e^{sx},\ \phi_{s,1}:=xe^{sx},\ \dots,\
\phi_{s,\nu-1}:=x^{\nu-1}e^{sx}\}
\end{equation}
while if $t=a+\sqrt{-1}b$ for $b>0$ is a complex root of multiplicity $\nu\ge2$,
then we must consider the family of functions:
\begin{equation}\label{eqn-1.d}
\begin{array}{l}
\{\phi_{t,0}:=e^{ax}\cos(bx),
\phi_{t,1}:=xe^{ax}\cos(bx),\dots,
\phi_{t,\nu-1}:=x^{\nu-1}e^{ax}\cos(bx),\\
\phantom{\{}\tilde\phi_{t,0}:=e^{ax}\sin(bx),
\tilde\phi_{t,1}:=xe^{ax}\sin(bx),\dots,
\tilde\phi_{t,\nu-1}:=x^{\nu-1}e^{ax}\sin(bx)\}\,.\vphantom{\vrule height 12pt}
\end{array}\end{equation}

\begin{theorem}\label{thm-1.4}
Assume that $r_+(P)>0$ and that $r_-(P)<0$.\begin{enumerate}
\item  If $s_1=r_+(P)$ and if the multiplicity of $s_1$ as a
root of $\mathcal{P}$ is larger than the corresponding multiplicity
of any complex root $t$ of $\mathcal{P}$ with $\Re(t)=s_1$, then
$\kappa_+[\sigma_{\Psi,P}]<\infty$; otherwise
$\kappa_+[\sigma_{\Psi,P}]=\infty$.
\item If $s_k=r_-(P)$ and if the multiplicity of $s_k$ as a
root of $\mathcal{P}$ is {{larger}} than the corresponding
multiplicity of any complex root $t$ of $\mathcal{P}$ with
$\Re(t)={s_k}$, then $\kappa_-[\sigma_{\Psi,P}]<\infty$; otherwise
$\kappa_-[\sigma_{\Psi,P}]=\infty$.
\end{enumerate}
\end{theorem}

\subsection{Indefinite signature inner products}
In Section~\ref{sect-7} we make some preliminary observations regarding
the situation if the inner product on $\mathbb{R}^n$ is non-degenerate with
indefinite signature.

We hope that these families of examples, which arise quite naturally from the
study of ordinary differential equations, help to shed light on questions of the
finiteness of the total first curvature.

\section{The proof of Theorem~\ref{thm-1.1}}\label{sect-2}
Assume all the roots of $\mathcal{P}$ are simple.
It then follows from the definition that
$$
||\sigma_P||^2=\sum_{i=1}^ke^{2s_ix}+
\sum_{j=1}^\ell e^{2a_jx}\,.
$$
Thus $||\sigma_P||^2$ tends to infinity as $x\rightarrow\infty$ if and only if some $s_i$
or some $a_j$ is positive or, equivalently, if $r_+(P)>0$. This implies that $\sigma_P$
is a proper map from $\HLP$ to $\mathbb{R}^n$ and that the length
is infinite. If $s_1>0$, then $\phi_1=e^{s_1x}$ is an injective map from
$\mathbb{R}$ to $\mathbb{R}$ and consequently
$\sigma_P$ is an embedding of $\mathbb{R}$ into $\mathbb{R}^n$.
If $a_1>0$, then $e^{a_1x}(\cos(b_1x),\sin(b_1x))$ is an injective map from
$\mathbb{R}$ to $\mathbb{R}^2$ and again we
may conclude that $\sigma_P$ is an embedding.
The analysis on $\HLN$ is similar if $r_-(P)<0$ and is therefore
omitted in the interests of brevity.\hfill\qed

\section{The proof of Theorem~\ref{thm-1.2}}\label{sect-3}

Throughout our proof, we will let $C_i=C_i(P)$ denote a generic positive
constant; we clear the notation after each case under consideration
and after the end of any given proof; thus $C_i$ can have
different meanings in different proofs or in different sections of the same proof.
We shall examine $\sigma_P$ on $\HLP$; the analysis
on $\HLN$ is similar and will therefore be omitted. We suppose $r_+>0$
or, equivalently, that $\max(s_1,a_1)>0$. We also assume that all the
roots of $\mathcal{P}$ are simple.
Suppose first that $s_1>a_1$ or that there are no complex roots. Let
$$
\epsilon:=\min_{\lambda\in\mathcal{R},\lambda\ne s_1}
(s_1-\Re(\lambda))=\min_{i>1,j\ge 1}(s_1-s_i,s_1-a_j)>0\,.$$
This measures the difference between the exponential growth rate of $\phi_1$
and the growth (or decay) rates of the functions
$\phi_i$ of Equation~(\ref{eqn-1.b}) for $i>1$
as $x\rightarrow\infty$. We have
\begin{equation}\label{eqn-3.a}
||\dot\sigma_P\wedge\ddot\sigma_P||^2=\sum_{i<j}
(\dot\phi_i\ddot\phi_j-\dot\phi_j\ddot\phi_i)^2\,.
\end{equation}
Consequently, we may estimate:
\begin{eqnarray}
&&||\dot\sigma_P\wedge\ddot\sigma_P||\le
C_1e^{(2s_1-\epsilon)x},\qquad
||\dot\sigma_P||^2\ge C_2e^{2s_1x}\text{ for }x\ge0,\nonumber\\
&&\frac{||\dot\sigma_P\wedge\ddot\sigma_P||}{||\dot\sigma_P||^2}\le
C_3e^{-\epsilon x}\text{ for }x\ge0\,.\label{eqn-3.b}
\end{eqnarray}
We integrate the estimate of Equation~(\ref{eqn-3.b}) to see $\kappa_+<\infty$.

Next suppose that $a_1>0$ and that $a_1\ge s_1$ (or that there are no real roots).
Then $e^{a_1x}$ is
the dominant term and we have
\begin{equation}\label{eqn-3.c}
||\dot\sigma_P||^2\le C_1e^{2a_1x}\,.
\end{equation}
The term $(\dot\phi_i\ddot\phi_j-\dot\phi_j\ddot\phi_i)^2$
in Equation~(\ref{eqn-3.a}) is maximized for $x\ge0$ when we take
$\phi_i=e^{a_1x}\cos(b_1x)$ and $\phi_j=e^{a_1x}\sin(b_1x)$.
We have:
\begin{eqnarray*}
&&\dot\phi_i=e^{a_1x}(a_1\cos(b_1x)-b_1\sin(b_1x))\\
&&\ddot\phi_i=e^{a_1x}\{(a_1^2-b_1^2)\cos(b_1x)-2a_1b_1\sin(b_1x)\}\\
&&\dot\phi_j=e^{a_1x}(a_1\sin(b_1x)+b_1\cos(b_1x))v_2,\\
&&\ddot\phi_j=e^{a_1x}\{(a_1^2-b_1^2)\sin(b_1x)+2a_1b_1\cos(b_1x)\},\\
&&\dot\phi_i^2+\dot\phi_j^2=(a_1^2+b_1^2)e^{2a_1x},\\
&&(\dot\phi_i\ddot\phi_j-\dot\phi_j\ddot\phi_i)^2
=b_1^2(a_1^2+b_1)^2e^{4a_1x}\,.
\end{eqnarray*}
Since $b_1\ne0$, we may estimate:
\begin{equation}\label{eqn-3.d}
||\dot\sigma_P\wedge\ddot\sigma_P||\ge C_2e^{2a_1x}\,.
\end{equation}
We use Equation~(\ref{eqn-3.c}) and Equation~(\ref{eqn-3.d}) to see
\begin{equation}\label{eqn-3.e}
\frac{||\dot\sigma_P\wedge\ddot\sigma_P||}{||\dot\sigma_P||^2}
\ge\frac{C_2}{C_1}>0\,.
\end{equation}
We integrate the uniform
estimate of Equation~(\ref{eqn-3.e}) to see $\kappa_+=\infty$.
\hfill\qed

\section{Examples}\label{sect-4} We now examine several specific cases. Since
the eigenvalues are to be simple, we can just specify
$\mathcal{P}$ or equivalently $\mathcal{R}$; the corresponding operator $P$ is then:
$$P=\mathcal{P}\left(\frac{d}{dx}\right)
=\prod_{\lambda\in\mathcal{R}}\left\{\frac{d}{dx}-\lambda\right\}\,.$$

\begin{example}\rm
Let $\mathcal{P}(\lambda)=\lambda^n-1$. The roots of $\mathcal{P}$
are the $n^{\operatorname{th}}$ roots of unity and all the roots
have multiplicity 1. Since $\mathcal{P}(1)=0$, $1$ is always a root.
\smallbreak\noindent{\bf Case I:} Suppose that $n$ is odd. Then $1$
is the only real root of $\mathcal{P}$. The remaining roots are all
complex. Thus $k=1$ and it follows that $\sigma_P$ is a proper
embedding of infinite length from $\HLP$ to $\mathbb{R}^n$.
If $\lambda^n=1$ and $\lambda\ne1$, then necessarily
$\Re(\lambda)<1$. It now follows that $\kappa_+$ is finite. There
exists a complex $n^{\operatorname{th}}$ root of unity with
$\Re(\lambda)<0$. Consequently, $\sigma_P$ is also a proper
embedding of infinite length from $\HLN$ to $\mathbb{R}^n$.
Since there are no real roots with $s_i<0$, we conclude $\kappa_-$
is infinite. \smallbreak\noindent{\bf Case II:} Suppose that $n$ is
even. Then $\pm1$ are the two real roots of $\mathcal{P}$. It now
follows that $\sigma_P$ is a proper embedding of infinite length
from $\HLP$ to $\mathbb{R}^n$ and from $\HLN$ to $\mathbb{R}^n$. If
$\lambda^n=1$ and $\lambda$ is not real, then $-1<\Re(\lambda)<1$.
Consequently, $\kappa_+$ and $\kappa_-$ are both finite.
\end{example}

\begin{example}
\rm Let $n\ge3$. Let $\{1,\dots,n-2,-1\pm\sqrt{-1}\}$ be the roots
of $\mathcal{P}$. Then $\sigma_P$ is a proper
embedding of infinite length from $\HLP$ to $\mathbb{R}^n$
and from $\HLN$ to $\mathbb{R}^n$,
 $\kappa_+$ is finite, and $\kappa_-$ is infinite.
\end{example}

\begin{example}
\rm Let $n\ge3$. Let $\{-1,\dots,2-n,1\pm\sqrt{-1}\}$ be the roots
of $\mathcal{P}$. Then $\sigma_P$ is a proper
embedding of infinite length from $\HLP$ to $\mathbb{R}^n$
and from $\HLN$ to $\mathbb{R}^n$,
$\kappa_+$ is infinite, and $\kappa_-$ is finite.
\end{example}

\begin{example}
\rm Let $n\ge2$. Let $\{1,\dots,n-1,-1\}$ be the roots of $\mathcal{P}$.
Then $\sigma_P$ is a proper
embedding of infinite length from $\HLP$ to $\mathbb{R}^n$
and from $\HLN$ to $\mathbb{R}^n$,
$\kappa_+$ is finite, and $\kappa_-$ is finite.
\end{example}

\begin{example}
\rm Let $n=2k\ge4$ be even. Let
$$\{1\pm\sqrt{-1},-1\pm\sqrt{-1},\dots,-(k-1)\pm\sqrt{-1}\}$$
be the roots of
$\mathcal{P}$.
Then $\sigma_P$ is a proper
embedding of infinite length from $\HLP$ to $\mathbb{R}^n$
and from $\HLN$ to $\mathbb{R}^n$,
$\kappa_+$ is infinite, and $\kappa_-$ is infinite.
\end{example}

\begin{example}
\rm Let $n=2k+1\ge5$ be odd. Let
$$\{0,1\pm\sqrt{-1},-1\pm\sqrt{-1},\dots,-(k-1)\pm\sqrt{-1}\}$$
be the roots of
$\mathcal{P}$. Then $\sigma_P$ is a proper
embedding of infinite length from $\HLP$ to $\mathbb{R}^n$
 and from $\HLN$ to $\mathbb{R}^n$,
$\kappa_+$ is infinite, and $\kappa_-$ is infinite.
\end{example}

\section{The proof of Theorem~\ref{thm-1.3}}\label{sect-5}

Let $\Phi=\{\phi_1,\dots,\phi_n\}$ be the standard basis for
$\mathcal{S}$ given in Equation~(\ref{eqn-1.b}) and let
$\Psi=\{\psi_1,\dots,\psi_n\}$ be any other basis for $\mathcal{S}$. Express
$$\psi_i=\Theta_i^j\phi_j$$ where we adopt
the {\it Einstein} convention and sum over repeated indices. We use $\Theta_i^j$
to make a linear change of basis on $\mathbb{R}^n$ and to regard
$\sigma_{\Psi,P}=\Theta\circ\sigma_P$; correspondingly, this defines a new
inner product $\langle\cdot,\cdot\rangle:=\Theta^*(\cdot,\cdot)$ on $\mathbb{R}^n$
so that
\begin{equation}\label{eqn-5.a}
||\dot\sigma_{\Psi,P}||=||\dot\sigma_P||_\Theta
\text{ and }
||\dot\sigma_{\Psi,P}\wedge\ddot\sigma_{\Psi,P}||=
||\dot\sigma_P\wedge\ddot\sigma_P||_\Theta\,.
\end{equation}
Any two norms on a finite dimensional real vector space are equivalent.
Thus
\begin{equation}\label{eqn-5.b}
C_1||v||\le ||v||_{\Theta}\le C_2||v||\,.
\end{equation}
The desired result now follows from Theorem~\ref{thm-1.1},
Theorem~\ref{thm-1.2}, Equation~(\ref{eqn-5.a}), and
Equation~(\ref{eqn-5.b}).\hfill\qed

\section{The proof of Theorem~\ref{thm-1.4}}\label{sect-6}
We will assume that $\Psi$ is the standard basis for $\mathcal{S}$
as the methods discussed in Section~\ref{sect-5} suffice to derive
the general result from this specific example.  We shall deal with
$\HLP$ as the situation for $\HLN$ is similar. The proof that
$r_+(P)>0$ implies $\sigma_P$ is a proper embedding of
$\HLP$ into $\mathbb{R}^n$ with infinite length is unchanged by any
questions of multiplicity since $e^{sx}$ or
$\{e^{ax}\cos(bx),e^{ax}\sin(bx)\}$ are still among the solutions of
$P$ for suitably chosen $s$ or $(a,b)$. We adopt the notation of
Equation~(\ref{eqn-1.c}) to define the functions
$\phi_{s,\mu}=x^{\mu}e^{sx}$ for $s\in\mathbb{R}$ and we adopt the
notation of Equation~(\ref{eqn-1.d}) to define the functions
$\phi_{t,\mu}=x^\mu e^{ax}\cos(bx)$ and $\tilde\phi_{t,\mu}=x^\mu
e^{ax}\sin(bx)$ for $t=a+b\sqrt{-1}$. We divide our discussion of
$\kappa_+$ into several cases: \smallbreak\noindent{\bf Case I:}
Suppose that $s_1>a_1$ and that $s_1$ is a real root of order $\nu$.
If $\nu=1$, the proof of Theorem~\ref{thm-1.2} extends to show
$\kappa_+$ is finite; the multiplicity of the other roots plays no
role as the exponential decay $e^{-\epsilon x}$ swamps any powers of
$x$. We suppose therefore that the multiplicity $\nu(s_1)>1$. We
will show that there exists $x_0$ so that:
\begin{eqnarray}\label{eqn-6.ax}
&&||\dot\sigma_P||^2\ge C_1x^{2\nu-2}e^{2s_1x}\text{ for }x\ge x_0,\\
&&||\dot\sigma_P\wedge\ddot\sigma_P||\le C_2x^{2\nu-4}e^{2s_1x}\text{ for }
x\ge x_0\,.
\label{eqn-6.bx}
\end{eqnarray}
It will then follow that
$$\frac{||\dot\sigma_P\wedge\ddot\sigma_P||}{||\dot\sigma_P||^2}\le C_3x^{-2}
\text{ for }x\ge x_0\,.$$
Since this is integrable on $\HLP$, we may conclude $\kappa_+$ is finite as desired.

We establish
Equation~(\ref{eqn-6.ax}) by noting that we have the following estimate:
\begin{eqnarray*}
||\dot\sigma_P||^2&=&\sum_{i=1}^n|\dot\phi_i|^2\ge|\dot\phi_{s_1,\nu-1}|^2
=\{s_1x^{\nu-1}+(\nu-1)x^{\nu-2}\}^2e^{2s_1x}\\
&\ge&s_1^2x^{2(\nu-1)}e^{2s_1x}\text{ for }x\text{ sufficiently large}\,.
\end{eqnarray*}
When dealing with $\HLP$, we may take $x_0=1$. However, when dealing with
$\HLN$, we must take $x_0<<0$ to ensure that the term $s_1x^{\nu-1}$
dominates the term $(\nu-1)x^{\nu-2}$ since these terms might have opposite signs
and cancellation could occur.

We may compute that:
\begin{equation}\label{eqn-6.c}
||\dot\sigma_P\wedge\ddot\sigma_P||^2=\sum_{i<j}
(\dot\phi_i\ddot\phi_j-\dot\phi_j\ddot\phi_i)^2\,.
\end{equation}
The assumption $s_1>a_1$ shows that the maximal term in this sum
occurs when $\phi_i=\phi_{s_1,\nu-1}$ and
$\phi_j=\phi_{s_1,\nu-2}$ and thus
$$||\dot\sigma_P\wedge\ddot\sigma_P||^2\le\frac{n(n-1)}2
\{\dot\phi_{s_1,\nu-1}\ddot\phi_{s_1,\nu-2}-\dot\phi_{s_1,\nu-2}\ddot\phi_{s_1,\nu-1}\}^2
\text{ for }x\ge x_0\,.$$ We have:
\begin{eqnarray*}
&&\dot\phi_{s_1,\nu-1}=(s_1x^{\nu-1}+(\nu-1)x^{\nu-2})e^{s_1x},\\
&&\ddot\phi_{s_1,\nu-1}=(s_1^2x^{\nu-1}+2s_1(\nu-1)x^{\nu-2}
+(\nu-1)(\nu-2)x^{\nu-3})e^{s_1x},\\
&&\dot\phi_{s_1,\nu-2}=(s_1x^{\nu-2}+(\nu-2)x^{\nu-3})e^{s_1x},\\
&&\ddot\phi_{s_1,\nu-2}=(s_1^2x^{\nu-2}+2s_1(\nu-2)x^{\nu-3}
+(\nu-2)(\nu-3)x^{\nu-4})e^{s_1x},\\
\end{eqnarray*}
Consequently:
\begin{eqnarray*}
&&\dot\phi_{s_1,\nu-1}\ddot\phi_{s_1,\nu-2}-\dot\phi_{s_1,\nu-2}\ddot\phi_{s_1,\nu-1}\\
&=&\{(s_1x^{\nu-1}+(\nu-1)x^{\nu-2})\\
&&\qquad\times
(s_1^2x^{\nu-2}+2s_1(\nu-2)x^{\nu-3}+(\nu-2)(\nu-3)x^{\nu-3})\}e^{2s_1x}\\
&-&\{(s_1x^{\nu-2}+(\nu-2)x^{\nu-3})\\
&&\qquad\times(s_1^2x^{\nu-1}+2s_1(\nu-1)x^{\nu-2}
+(\nu-1)(\nu-2)x^{\nu-3})e^{2s_1x}
\end{eqnarray*}
The leading terms cancel:
$$\{(s_1x^{\nu-1}s_1^2x^{\nu-2})-(s_1x^{\nu-2}s_1^2x^{\nu-1})\}e^{2s_1x}=0\,.$$
The remaining terms are $O(x^{2\nu-4}e^{2s_1x})$ as desired;
Equation~(\ref{eqn-6.bx}) now follows. This shows $\kappa_+$ is finite if
$s_1>a_1$.

\medbreak\noindent{\bf Case II:} Suppose $a_1>s_1$. Choose the
complex root $t_1=a_1+b_1\sqrt{-1}$ to have maximal multiplicity
$\nu$ among all the complex roots $t\in\mathcal{R}$ with
$\Re(t)=a_1$. The dominant term in Equation~(\ref{eqn-6.c}) occurs
when $\phi_i=\phi_{t_1,\nu-1}$ and
$\phi_j=\tilde\phi_{t_1,\nu-1}$. Differentiating powers of
$x$ lowers the order in $x$ and give rise to lower order terms. Thus
we may ignore these derivatives and use the computations performed
in Section~\ref{sect-3} to see:
\begin{eqnarray*}
&&C_1x^{2\nu-2}e^{2a_1x}\le||\dot\sigma_P||^2\le
C_2x^{2\nu-2}e^{2a_1x}
\text{ for }x\ge x_0,\\
&&(\dot\phi_i\ddot\phi_j-\dot\phi_j\ddot\phi_i)^2\ge
C_3x^{4(\nu-1)}e^{4a_1x} \text{ for }x\ge
x_0\,.
\end{eqnarray*}
We may now conclude that $\kappa_+=\infty$.

\medbreak\noindent{\bf Case III:} The difficulty comes when $a_1=s_1$.
If $t_1$ is a complex root of multiplicity at least
as great as
 the multiplicity of $s_1$, the $\{\phi_{t_1,\nu-1},\tilde\phi_{t_1,\nu-1}\}$
 terms dominate the
 computation and the argument given above in Case II
 implies $\kappa_+$ is infinite. On the other hand, if all the complex roots
 with $\Re(\lambda)=s_1$ have
 multiplicity less than the multiplicity of $s_1$, then the
 $\phi_{s_1,\nu-1}$ 
 terms dominate the computation and the argument given above in Case I
 shows that $\kappa_+$ is finite.
 \hfill\qed

\bigbreak We conclude this section with an example where the
multiplicity plays
 a crucial role and where our previous results are not applicable.
 \begin{example}\rm
 Let $P(y)=y^{(n)}$ for $n\ge2$. Then $\mathcal{R}=\{0\}$ and $0$ is
 a root of multiplicity $n$. We have
 $\mathcal{S}=\operatorname{Span}\{\phi_1:=1,\phi_2:=x,...,\phi_{n}:=x^{n-1}\}$. Since $x\in\mathcal{S}$,
 $\sigma_P$ is a proper map of infinite length on $\HLP$ and on $\HLN$. We
 have:
 \begin{eqnarray*}
&&|\dot\sigma_P|^2\ge C_1x^{2n-2},\text{ and }\\
&& \sum_{i<j}(\ddot\phi_i\dot\phi_j-\ddot\phi_j\dot\phi_i)^2
=\sum_{i<j}((i-1)(i-2)(j-1)-(j-1)(j-2)(i-1))^2)x^{2(i+j-3)}\\
&&\qquad\le C_2x^{2(2n-4)}\,.
 \end{eqnarray*}
 Consequently
 $|\kappa|\le C_3\frac{x^{2n-4}}{x^{2n-2}}\text{ for }|x|\ge1$.
 This is integrable so $\kappa_+<\infty$ and $\kappa_-<\infty$.
 \end{example}\black

\section{Higher signature geometry}\label{sect-7}
According to modern physics, we do not live in a Riemannian universe, but
rather in a Lorentzian universe or even in the higher signature setting according
to some string theories.
Indefinite signatures also appears in supergravity theory (see, for example,
Abounasr, Belhaj, Rasmussen, and Saidi \cite{ABRS06}). Walker
metrics have indefinite signatures and are important as well (see, for example,
the discussion in Davidov, Diaz-Ramos, Garcia-Rio, Matsushita, Muskarov,
and Vazquez-Lorenzo \cite{D08}).

Thus it is natural to consider the setting where
$\langle\cdot,\cdot\rangle$ is a non-degenerate symmetric bilinear form of
signature $(p,q)$ on the underlying vector space $\mathbb{R}^n$. The geometry
here is very different and warrants investigation in its own right; we shall content
ourselves with just a preliminary remarks in this current paper and postpone until
a subsequent paper a more detailed investigation.

The existence of null directions (i.e. vectors $0\ne v$
where $\langle v,v\rangle=0$) causes some additional technical complications
that are at the heart of the matter; thus, for example, it is not possible to use
Equation~(\ref{eqn-1.a}) to define the first curvature $\kappa$ where the tangent
vector $\dot\kappa$ is null.

For the remainder of this section, we shall work in signature $(1,1)$.  We consider
the operator
$P(y):=y^{(2)}-y$ so that $\mathcal{S}=\operatorname{Span}\{e^x,e^{-x}\}$
and $\sigma_P=(e^{x},e^{-x})$. This curve is a proper embedding.
We let $\langle\cdot,\cdot\rangle$ be
the inner product on $\mathbb{R}^2$ defined by the symmetric matrix:
$$G:=\left(\begin{array}{ll}\alpha&\beta\\\beta&\gamma\end{array}\right)\,.$$
This means that:
$$\langle (w_1,w_2),(z_1,z_2)\rangle=\alpha w_1z_1
+\beta(w_1z_2+w_2z_1)+\gamma z_2w_2\,.$$
We suppose $\det(G)<0$ to ensure $G$ has signature $(1,1)$.
We have:
\begin{eqnarray*}
&&\dot\sigma_P=(e^{x},-e^{-x}),\\
&&\langle\dot\sigma_P,\dot\sigma_P\rangle=\alpha e^{2x}
+\gamma e^{-2x}-2\beta,\\
&&\ddot\sigma_P=(e^{x},e^{-x}),\\
&&\dot\sigma_P\wedge\ddot\sigma_P=2
\{(1,0)\wedge(0,1)\},\\
&&\langle\dot\sigma_P\wedge\ddot\sigma_P,
\dot\sigma_P\wedge\ddot\sigma_P\rangle=4
(\alpha\gamma-\beta^2)
\end{eqnarray*}
Since the induced metric on $\Lambda^2(\mathbb{R})$ is {\it negative} definite,
$\alpha\gamma-\beta^2$ is negative and it is natural to take the absolute value
when extracting the square root:
$$
\kappa ds:=\frac{2|\alpha\gamma-\beta^2|^{1/2}}
{|\alpha e^{2x}+\gamma e^{-2x}-2\beta |}
dx\,.
$$

\begin{example}\rm
Let $\alpha=\gamma=0$ and $\beta=-1$.
Then $\sigma_P$ is spacelike, of
infinite length on $\HLP$ and on $\HLN$, $\kappa_+=\infty$, and $\kappa_-=\infty$.
\end{example}

\begin{example}\rm
Let $\alpha=\gamma=\epsilon>0$ and $\beta=-1$ where $\epsilon$ is small.
Then $\sigma_P$ is spacelike, of infinite length on $\HLP$ and on $\HLN$,
$\kappa_+<\infty$, and $\kappa_-<\infty$.\end{example}

\begin{example}\rm
Let $\alpha=\epsilon$, $\beta=-1$, and $\gamma=0$ where $\epsilon$
is small. Then $\sigma_P$ is spacelike, of infinite length
on $\HLP$ and on $\HLN$, $\kappa_+<\infty$, and
$\kappa_-=\infty$.
\end{example}

\begin{example}\rm
Let $\alpha=0$, $\beta=-1$, and $\gamma=\epsilon$ where $\epsilon$
is small. Then $\sigma_P$ is spacelike, of infinite length
on $\HLP$ and on $\HLN$, $\kappa_+=\infty$, and $\kappa_-<\infty$.
\end{example}

\section*{Acknowledgements}
This work was partially supported by the National Research
Foundation of Korea (NRF) grant funded by the Korea government
(MEST) (2012-0005282) and by project MTM2009-07756 (Spain).
\end{document}